\documentclass[12pt]{article}

\usepackage{booktabs}
\usepackage{amsmath}
\usepackage{amsfonts}
\usepackage{amssymb}
\usepackage{amsthm}
\usepackage{url}

\usepackage{amsthm}

\newtheorem{theorem}{Theorem}
\newtheorem{corollary}{Corollary}
\newtheorem{lemma}{Lemma}

\theoremstyle{definition}
\newtheorem{definition}{Definition}

\newcommand{\pgl}{\operatorname{PGL}}
\newcommand{\pgammal}{\operatorname{P}\!\Gamma\!\operatorname{L}}
\newcommand{\pg}{\operatorname{PG}}

\newcommand{\gf}{\operatorname{GF}}
\newcommand{\aut}{\operatorname{Aut}}
\newcommand{\con}{\operatorname{Conj}}
\newcommand{\B}{{\mathcal B}}

\begin{document}

\title{Update: Some new results on lower bounds on $(n,r)$-arcs in $\pg(2,q)$ for $q\le 31$}
\author{Michael Braun (\url{michael.braun@h-da.de})\\
Faculty of Computer Science\\
University of Applied Sciences, Darmstadt, Germany}

\maketitle

\begin{abstract}
An $(n,r)$-arc in $\pg(2,q)$ is a set $\B$ of points in $\pg(2,q)$ such that each line in $\pg(2,q)$ contains at most $r$ elements of $\B$ and such that there is at least one line containing exactly $r$ elements of $\B$. The value $m_r(2,q)$ denotes the maximal number $n$ of points in the projective geometry $\pg(2,q)$ for which an $(n,r)$-arc exists. By explicitly constructing $(n,r)$-arcs using prescribed automorphisms and integer linear programming we obtain some improved lower bounds for $m_r(2,q)$: $m_{10}(2,16)\ge 144$, $m_3(2,25)\ge 39$, $m_{18}(2,25)\ge 418$, $m_9(2,27)\ge 201$, $m_{14}(2,29)\ge 364$, $m_{25}(2,29)\ge 697$, $m_{25}(2,31)\ge 734$. Furthermore, we show by systematically excluding possible automorphisms that putative $(44,5)$-arcs, $(90,9)$-arcs in $\pg(2,11)$, and $(39,4)$-arcs in $\pg(2,13)$---in case of their existence---are rigid, i.e. they all would only admit the trivial automorphism group of order $1$. In addition, putative $(50,5)$-arcs, $(65,6)$-arcs, $(119,10)$-arcs, $(133,11)$-arcs, and $(146,12)$-arcs in $\pg(2,13)$ would be rigid or would admit a unique automorphism group (up to conjugation) of order $2$.
\end{abstract}

\section{Introduction}

\begin{definition}
An \emph{$(n,r)$-arc} in $\pg(2,q)$ is a set $\B$ of points in $\pg(2,q)$ such that each line in $\pg(2,q)$ contains at most $r$ elements of $\B$ and such that there is at least one line containing exactly $r$ elements of $\B$. 
\end{definition}

It is well-known (e.g. see \cite{BH05}) that $(n,r)$-arcs in $\pg(2,q)$ are closely related to error-correcting linear codes: The $n$ points of an $(n,r)$-arc in $\pg(2,q)$ define the columns of a $3\times n$ generator matrix of linear $[n,3,n-r]_q$ code, which is code of length $n$, dimension $3$, and minimum distance $n-r$ with respect to the Hamming metric. The linear code is projective since the columns of any generator matrix are pairwise linearly independent.

\begin{definition}
Let $m_r(2,q)$ denote the maximum number $n$ for which an $(n,r)$-arc in $\pg(2,q)$ exists.
\end{definition}

A major goal in studying $(n,r)$-arcs in $\pg(2,q)$ is the determination of $m_r(2,q)$.

In general it is hard to determine the exact value of $m_r(2,q)$ and in most cases instead of the exact value only a lower and an upper bound for $m_r(2,q)$ are known. An explicit construction of an $(n,r)$-arc in $\pg(2,q)$ yields a lower bound $m_r(2,q)\ge n$.

The values $m_r(2,q)$ with $q\le 9$ are exactly determined (see \cite{BH05}). For $m_r(2,q)$ with $11\le q\le 19$ we refer to \cite{Bal} whereas a table for $23\le q\le 31$ can compiled from several sources \cite{Bra19,CH,Das08,DM11,DM13,DM17a,DM17b,DM17c,DM18,Das19,DM20,DGM+11,HS01}. An recent overview with tables on all values $q\le31$ can be found in \cite{Bra19}. 

In this article, we extend the results from \cite{Bra19} of lower bounds on $m_r(2,q)$ and give some improvements listed in Table~\ref{tab:results}.

\begin{table}[!htbp]
\centering
\caption{Improved lower bounds on $m_r(2,q)$}\label{tab:results}
\begin{tabular}{llrr}
\toprule
$q$ & $r$ & old bound & new bound\\
\midrule
$16$ & $10$ & $142$ & $144$\\
$25$ & $3$ & $38$ & $39$\\ 
$25$ & $18$ & $416$ & $418$\\
$27$ & $9$ & $198$ & $201$\\
$29$ & $14$ & $361$ & $364$\\
$29$ & $25$ & $695$ & $697$\\
$31$ & $25$ & $733$ & $734$\\
\bottomrule
\end{tabular}
\end{table}

Furthermore, as a second result in this work, we show that the putative $(n,r)$-arcs in $\pg(2,q)$ for $q\in\{11,13\}$ for the open gaps between lower and upper bound on $m_r(2,q)$---in case of existence---are rigid or only admit a unique automorphism of order $2$.

\section{Construction by integer linear programming}

We use the construction of $(n,r)$-arcs in $\pg(2,q)$ with prescribed groups of automorphisms using integer linear programming described in \cite{BKW05,Bra19}:

In the following, let 
\[
{\gf(q)^n\brack k}
\]
denote the set of $k$-dimensional subspaces of $\gf(q)^n$, which is called the \emph{Grassmannian}. Its cardinality is given by the \emph{Gaussian number}, also called $q$-Binomial coefficient:
\[
{n\brack k}_q=\left|{\gf(q)^n\brack k}\right| 
=\prod_{i=0}^{k-1}\frac{q^n-q^i}{q^k-q^i}.
\]

In terms of vector spaces, an $(n,r)$-arc in $\pg(2,q)$ corresponds to a set $\B\subseteq{\gf(q)^3\brack 1}$ such that for all $H\in{\gf(q)^3\brack 2}$ holds:
\[
|\{P\in\B\mid H\supseteq P\}|\le r.
\]

If ${\gf(q)^3\brack 1}=\{P_1,\ldots,P_{q^2+q+1}\}$ and ${\gf(q)^3\brack 2}=\{H_1,\ldots,H_{q^2+q+1}\}$, where ${3\brack 1}_q={3\brack 2}_q=q^2+q+1$, we define the $(q^2+q+1)\times (q^2+q+1)$ incidence matrix 
\[
A(q)=(a_{ij})
\]
with entries
\[
a_{ij}:=\begin{cases}
1&\text{if $H_i\supseteq P_j$,}\\
0&\text{otherwise.}
\end{cases}
\]

\begin{lemma}
If $u=(1,\ldots,1)^T$ denote the all-one vector any binary column vector $x$ satisfying
\[
A(q)\cdot x\le r\cdot u
\]
is equivalent to a $(u^T\cdot x,r)$-arc in $\pg(2,q)$.
\end{lemma}

\begin{corollary}
The determination of $m_r(2,q)$ corresponds to the following integer linear programming problem:
\[
m_r(2,q)=\max_{x\in\{0,1\}^{q^2+q+1}}\{u^T\cdot x\mid A(q)\cdot x\le r\cdot u\}.
\]
\end{corollary}

The incidence preserving bijections (automorphisms) of our ambient space for $(n,r)$ arcs---the projective geometry $\pg(2,q)$---are defined by the projective semi-linear group $\pgammal(3,q)$ (see \cite{Bae52}). It acts transitively on the Grassmannian ${\gf(q)^3\brack k}$. 

Hence, any subgroup $G\le\pgammal(3,q)$ partitions the Grassmannian into $G$-orbits. If $\alpha\in\pgammal(3,q)$ and $S\in{\gf(q)^3\brack k}$ we denote by 
\[
\alpha S:=\{\alpha x\mid x\in S\}
\]
the transformed subspace and by 
\[
G(S):=\{\alpha S\mid \alpha\in G\}
\]
the $G$-orbit of $S$. The set of all $G$-orbits will be written as 
\[
G\backslash\!\!\backslash{\gf(q)^3\brack k}.
\]

\begin{definition}
An $(n,r)$-arc $\B$ in $\pg(2,q)$ admits a subgroup $G\le\pgammal(3,q)$ as a \emph{group of automorphisms} if and only if $\B$ consists of $G$-orbits on ${\gf(q)^3\brack 1}$. The maximal group of automorphisms of $\B$ is called the \emph{automorphism group} of $\B$ and abbreviated by 
\[
\aut(\B).
\]
\end{definition}

\begin{definition}
Let $m^G_r(2,q)$ denote the maximal size $n$ of an $(n,r)$-arc in $\pg(2,q)$ admitting $G\le\pgammal(3,q)$ as a group of automorphisms. 
\end{definition}

\begin{corollary}
For any $G\le\pgammal(3,q)$ we get a lower bound 
\[
m^G_r(2,q)\le m_r(2,q).
\]
In particular, for the trivial group $G=\{1\}$ we have 
\[
m^{\{1\}}_r(2,q)= m_r(2,q).
\]
\end{corollary}

If $\{P_1,\ldots,P_\ell\}$ denotes a set of representatives of the orbits $G\backslash\!\!\backslash{\gf(q)^3\brack 1}$ and $\{H_1\ldots,H_\ell\}$ a transversal of the orbits $G\backslash\!\!\backslash{\gf(q)^3\brack 2}$ for any $G\le\pgammal(3,q)$ we define the $G$-incidence matrix $A(G)=(a_{ij})$ with
\[
a_{ij}:=|\{P\in G(P_j)\mid H_i\supseteq P\}|.
\]
Furthermore, by $w(G)=(w_1,\ldots,w_\ell)^T$ we denote the vector of the lengths of $G$-orbits on ${\gf(q)^3\brack 1}$, i.e. 
\[
w_j:=|G(P_j)|.
\] 
Note that the number of orbits of $G$ on the set of points and hyperplanes is equal
\[
\ell=\left|G\backslash\!\!\backslash{\gf(q)^3\brack 1}\right|=
\left|G\backslash\!\!\backslash{\gf(q)^3\brack 2}\right|\le q^2+q+1.
\]

\begin{theorem}\label{thm:ilp}
Any binary vector $x$ of length $\ell=|G\backslash\!\!\backslash{\gf(q)^3\brack 1}|$ with  
\[
A(G)\cdot x\le r\cdot u
\]
corresponds to a $(w(G)^T\cdot x,r)$-arc in $\pg(2,q)$ admitting $G\le \pgammal(3,q)$ as a group of automorphisms. In addition, we obtain the following integer linear programming:
\[
m^G_r(2,q)=\max_{x\in\{0,1\}^\ell}\{w(G)^T\cdot x\mid A(G)\cdot x\le r\cdot u\}.
\]
\end{theorem}

\section{Constructed arcs}

We list $(n,r)$-arcs in $\pg(2,q)$ in the appendix constructed with the proposed approach and with improved size. All $(n,r)$-arcs in $\pg(2,q)$ were computed with Gurobi (see \cite{Gur}) as ILP solver. 

Elements of the prime field $\gf(p)$ are represented by integers $0\le a < p$ where elements are added and multiplied modulo $p$. In extension fields $\gf(p^e)$ the elements $\sum_{i=0}^{e-1}a_ix^i$ are given by the numbers $\sum_{i=0}^{e-1}a_ip^i$ where elements are added and multiplied modulo a given irreducible polynomial $f(x)\in\gf(p)[x]$ of degree $e$. For the finite fields $\gf(16)$, $\gf(25)$, and $\gf(27)$ we use the irreducible polynomials:
\begin{align*}
&x^4+x^3+1\in\gf(2)[x],\\
&x^2+x+2\in\gf(5)[x],\\
&x^3+2x+1\in\gf(3)[x].
\end{align*}

As group we used subgroups projective linear groups $G\le\pgl(3,q)\le\pgammal(3,q)$. For two cases we used the symmetric group $S_3$ generated by two matrices
\[
\begin{pmatrix}
0 & 1 & 0\\ 
0 & 0 & 1\\
1 & 0 & 0\\
\end{pmatrix}\text{ and }
\begin{pmatrix}
0 & 1 & 0\\ 
1 & 0 & 0\\
0 & 0 & 1\\
\end{pmatrix}.
\]
For the remaining cases we used random cyclic groups generated by an invertible $3\times 3$ matrix.

\section{Excluding automorphisms}

An open gap is an entry in the tables of $m_r(2,q)$ for which upper and lower bound differ:
\[
\ell\le m_r(2,q)\le u\quad\text{where}\quad \ell< u.
\]
In that case the question is whether an $(\ell+1,r)$-arc in $\pg(2,q)$ exits or not. We call such an arc a \emph{putative} arc in $\pg(2,q)$. In \cite{Bra} it was shown that for the gap $100\le m_{10}(2,11)\le101$ a putative $(101,10)$-arc in $\pg(2,11)$ admits---in case of its existence---only the trivial automorphism group of order $1$. 

In this paper we consider the remaining gaps for $q=11$ and $q=13$ which are given in Table~\ref{tab:open}.

\begin{table}[!htbp]
\centering
\caption{Open cases of $m_r(2,q)$ for $q=11,13$}\label{tab:open}
\begin{tabular}{rcl}
\toprule
$m_5(2,11)$&$\in$&$\{43,44,45\}$\\
$m_9(2,11)$&$\in$&$\{89,90\}$\\
$m_4(2,13)$&$\in$&$\{38,39,40\}$\\
$m_5(2,13)$&$\in$&$\{49,50,51,52,53\}$\\
$m_6(2,13)$&$\in$&$\{64,65,66\}$\\
$m_{10}(2,13)$&$\in$&$\{118,119\}$\\
$m_{11}(2,13)$&$\in$&$\{132,133\}$\\
$m_{12}(2,13)$&$\in$&$\{145,146,147\}$\\
\bottomrule
\end{tabular}
\end{table}

\begin{definition}
Let $\B,\B'$ be an $(n,r)$-arcs in $\pg(2,q)$ The two sets $\B$ and $\B'$ are defined to be isomorphic if and only if there exists $\alpha\in\pgammal(3,q)$ such that
\[
\alpha\B:=\{\alpha P\mid P\in\B\}=\B'
\]
The set of all arcs that are \emph{isomorphic} to $\B$ is denoted by
\[
\pgammal(3,q)(\B):=\{\alpha\B\mid \alpha\in\pgammal(3,q)\}.
\]
\end{definition}

Note that due to the incidence preserving property of $\pgammal(3,q)$ isomorphic arcs have the same parameters.

The following lemma is well-known from the theory of group actions (cf.~\cite{Ker99}) and states that the automorphism groups of isomorphic objects are conjugated.

\begin{lemma}
Let $\B$ be an $(n,r)$-arc in $\pg(2,q)$ and let $\alpha\in\pgammal(3,q)$. Then we obtain:
\[
\aut(\alpha\B)=\alpha\aut(\B)\alpha^{-1}=\{\alpha\beta\alpha^{-1}\mid \beta\in\aut(\B)\}.
\]
\end{lemma}

If $\B$ in an $(n,r)$-arc in $\pg(2,q)$ with $G\le\aut(\B)$ then any isomorphic arc $\B'=\alpha \B$ for $\alpha\in\pgammal(3,q)$
admits the conjugated group $G'=\alpha G\alpha^{-1}$ satisfies
\[
G'=\alpha G\alpha^{-1}\le\alpha\aut(G)\alpha^{-1}=\aut(\alpha\B)=\aut(\B'),
\]
which means that the conjugated group $G'$ also occurs as a group of automorphisms of $\B'$.

As a consequence, when aiming for $(n,r)$-arcs in $\pg(2,q)$ with prescribed groups of automorphisms it is sufficient to consider representatives of conjugacy classes of subgroups of $\pgl(3,q)$ as possible candidates for potential groups to be prescribed.

Furthermore, any $(n,r)$-arc $\B$ in $\pg(2,q)$ with $\{1\}<G\le\aut(\B)$ also admits all cyclic subgroups $C\le G$ as groups of automorphisms. 

\begin{corollary}
If we can show for all representatives $C$ of conjugacy classes of nontrivial cyclic subgroups of $\pgammal(3,q)$ that no $(n,r)$-arc in $\pg(2,q)$ exists with $C$ as group as automorphisms, either the automorphism group of such arcs are trivial or arcs with that set of parameters do not exist.
\end{corollary}

In case of a prime field $\gf(q)$ the projective semi-linear group is exactly the projective linear group
\[
\pgammal(3,q)=\pgl(3,q).
\]
In the following, a transversal of conjugacy classes of cyclic subgroups of $\pgl(3,q)$ will be abbreviated by 
\[
\con(q).
\]
Its cardinality is given by (see \cite{Mac81}):
\[
|\con(q)|=\begin{cases}
q^2+q+2&\text{if $3$ divides $q-1$,}\\
q^2+q&\text{otherwise.}
\end{cases}
\]

\begin{lemma}\label{lem:q11}
Let $q$ be a prime. If
\[
m^C_r(2,q)<n\quad \forall C\in\con(q)\setminus\{\{1\}\}
\] 
one of the following conditions holds:
\begin{enumerate}
\item $m_r(2,q)<n$.
\item $(n,r)$-arcs $\B$ in $\pg(2,q)$ exist where $\aut(\B)=\{1\}$.
\end{enumerate}
\end{lemma}

We now apply this corollary to the parameters $(q,n,r)=(11,44,5)$, $(q,n,r)=(11,90,9)$, and $(q,n,r)=(13,50,5)$. There are $|\con(11)|=132$ conjugacy classes of cyclic subgroups of $\pgl(3,11)$ and $|\con(13)|=184$ classes in $\pgl(3,13)$. We compute the representatives using GAP \cite{GAP}. By solving the integer linear programming $m^C_r(2,q)$ according to Theorem~\ref{thm:ilp} for all $C\in\con(q)\setminus\{1\}$ using Gurobi \cite{Gur} we obtain with a runtime less than 3 hours on a 1.2 GHz Intel Core m3 processor the following result:

\begin{theorem}
In case of their existence the automorphism groups of $(44,5)$-arcs, $(90,9)$-arcs in $\pg(2,11)$, and $(50,5)$-arcs in $\pg(2,13)$ would be trivial of order $1$.
\end{theorem}

For the remaining parameters $r\in\{5,6,10,11,12\}$ for $q=13$ we apply a slightly adapted version of Lemma~\ref{lem:q11} since for these open cases exactly one cyclic subgroup 
\[
C_0:=\langle
\begin{pmatrix}
0&1&0\\
1&0&0\\
0&0&12
\end{pmatrix}
\rangle
\]
of order $2$ could not directly be excluded to be a group of automorphisms of putative arcs since we cancelled the ILP solver for $m_r^{C_0}(2,13)$ after a 5000 seconds (for each value $r$).

But it is obvious to conject that this group can also be excluded if we spend more running time on the ILP solver.

\begin{lemma}
Let $q$ be a prime. Let $C_0\in\con(q)$. If
\[
m^C_r(2,q)<n\quad \forall C\in\con(q)\setminus\{\{1\},C_0\}
\] 
one of the following conditions holds:
\begin{enumerate}
\item $m_r(2,q)<n$.
\item $(n,r)$-arcs $\B$ in $\pg(2,q)$ exist where either $\aut(\B)=\{1\}$ or $\aut(\B)$ is conjugated to $C_0$.
\end{enumerate}
\end{lemma}

Finaly, we get

\begin{theorem}
In case of their existence the automorphism groups of $(50,5)$-arcs, $(65,6)$-arcs, $(119,10)$-arcs, $(133,11)$-arcs, and $(146,12)$-arcs in $\pg(2,13)$ would either be trivial of order $1$ or would be the following cyclic subgroup of order $2$ (up to conjugation):
\[
C_0:=\langle
\begin{pmatrix}
0&1&0\\
1&0&0\\
0&0&12
\end{pmatrix}.
\]
\end{theorem}

\section*{Appendix}

\subsection{$m_{10}(2,16)\ge 144$}

group: symmetric group $S_3$

\noindent
$(0,0,1)$, $(1,0,0)$, $(0,1,0)$, $(0,1,4)$, $(1,0,4)$, $(1,0,6)$, $(0,1,6)$, 
$(1,4,0)$, $(1,6,0)$, $(0,1,5)$, $(1,0,5)$, $(1,0,15)$, $(0,1,15)$, $(1,5,0)$, 
$(1,15,0)$, $(1,1,2)$, $(1,12,12)$, $(1,2,1)$, $(1,1,3)$, $(1,8,8)$, $(1,3,1)$, 
$(1,1,5)$, $(1,15,15)$, $(1,5,1)$, $(1,1,7)$, $(1,14,14)$, $(1,7,1)$, $(1,1,8)$, 
$(1,3,3)$, $(1,8,1)$, $(1,1,9)$, $(1,13,13)$, $(1,9,1)$, $(1,1,10)$, $(1,11,11)$, 
$(1,10,1)$, $(1,2,3)$, $(1,12,13)$, $(1,8,9)$, $(1,3,2)$, $(1,13,12)$, $(1,9,8)$, 
$(1,2,4)$, $(1,12,2)$, $(1,6,12)$, $(1,4,2)$, $(1,2,12)$, $(1,12,6)$, $(1,2,7)$, 
$(1,12,15)$, $(1,14,5)$, $(1,7,2)$, $(1,15,12)$, $(1,5,14)$, $(1,2,8)$, $(1,12,4)$, 
$(1,3,6)$, $(1,8,2)$, $(1,4,12)$, $(1,6,3)$, $(1,2,10)$, $(1,12,5)$, $(1,11,15)$, 
$(1,10,2)$, $(1,5,12)$, $(1,15,11)$, $(1,2,14)$, $(1,12,7)$, $(1,7,14)$, $(1,14,2)$, 
$(1,7,12)$, $(1,14,7)$, $(1,2,15)$, $(1,12,11)$, $(1,5,10)$, $(1,15,2)$, $(1,11,12)$, 
$(1,10,5)$, $(1,3,4)$, $(1,8,11)$, $(1,6,10)$, $(1,4,3)$, $(1,11,8)$, $(1,10,6)$, 
$(1,3,7)$, $(1,8,10)$, $(1,14,11)$, $(1,7,3)$, $(1,10,8)$, $(1,11,14)$, $(1,3,9)$, 
$(1,8,7)$, $(1,13,14)$, $(1,9,3)$, $(1,7,8)$, $(1,14,13)$, $(1,3,10)$, $(1,8,6)$, 
$(1,11,4)$, $(1,10,3)$, $(1,6,8)$, $(1,4,11)$, $(1,3,11)$, $(1,8,14)$, $(1,10,7)$, 
$(1,11,3)$, $(1,14,8)$, $(1,7,10)$, $(1,4,6)$, $(1,6,13)$, $(1,4,9)$, $(1,6,4)$, 
$(1,13,6)$, $(1,9,4)$, $(1,4,14)$, $(1,6,15)$, $(1,7,5)$, $(1,14,4)$, $(1,15,6)$, 
$(1,5,7)$, $(1,4,15)$, $(1,6,9)$, $(1,5,13)$, $(1,15,4)$, $(1,9,6)$, $(1,13,5)$, 
$(1,5,9)$, $(1,15,10)$, $(1,13,11)$, $(1,9,5)$, $(1,10,15)$, $(1,11,13)$, $(1,5,11)$, 
$(1,15,13)$, $(1,10,9)$, $(1,11,5)$, $(1,13,15)$, $(1,9,10)$, $(1,7,13)$, $(1,14,9)$, 
$(1,9,13)$, $(1,13,7)$, $(1,9,14)$, $(1,13,9)$

\subsection{$m_{3}(2,25)\ge 39$}

group: cyclic; generated by 
\[
\begin{pmatrix}
0&10&13 \\
13&3&5 \\
18&18&4
\end{pmatrix}
\]

\noindent
$(0,1,12)$, $(1,5,2)$, $(1,6,14)$, $(1,0,6)$, $(1,4,7)$, $(1,10,21)$, $(1,0,13)$, 
$(1,5,6)$, $(1,12,14)$, $(1,0,18)$, $(1,14,0)$, $(1,13,8)$, $(1,2,21)$, $(1,13,4)$, 
$(1,18,19)$, $(1,2,22)$, $(1,6,12)$, $(1,14,13)$, $(1,3,2)$, $(1,22,10)$, $(1,12,12)$, 
$(1,3,10)$, $(1,7,23)$, $(1,23,0)$, $(1,3,20)$, $(1,11,13)$, $(1,4,10)$, $(1,6,3)$, 
$(1,10,20)$, $(1,20,21)$, $(1,12,15)$, $(1,13,12)$, $(1,14,15)$, $(1,24,19)$, $(1,23,20)$, 
$(1,19,0)$, $(1,19,7)$, $(1,19,2)$, $(1,23,16)$

\subsection{$m_{18}(2,25)\ge 418$}

group: cyclic; generated by 
\[
\begin{pmatrix}
19& 24& 11 \\
21& 3& 18 \\
22& 8& 7 
\end{pmatrix}
\]

\noindent
$(0,1,0)$, $(1,20,9)$, $(1,21,5)$, $(1,5,12)$, $(0,1,1)$, $(1,10,2)$, $(1,10,12)$, 
$(1,15,9)$, $(0,1,3)$, $(1,11,12)$, $(1,6,9)$, $(1,1,10)$, $(0,1,4)$, $(1,8,16)$, 
$(1,0,13)$, $(1,8,19)$, $(0,1,5)$, $(1,18,23)$, $(1,18,21)$, $(1,2,24)$, $(0,1,7)$, 
$(1,13,7)$, $(1,7,3)$, $(1,3,8)$, $(0,1,9)$, $(1,9,1)$, $(1,19,15)$, $(1,19,20)$, 
$(0,1,10)$, $(1,22,4)$, $(1,17,2)$, $(0,1,14)$, $(0,1,11)$, $(1,2,15)$, $(1,9,16)$, 
$(1,9,3)$, $(0,1,12)$, $(1,12,22)$, $(1,14,18)$, $(1,21,4)$, $(0,1,15)$, $(1,0,20)$, 
$(1,5,10)$, $(1,24,6)$, $(0,1,16)$, $(1,1,5)$, $(1,13,24)$, $(1,7,5)$, $(0,1,17)$, 
$(1,24,24)$, $(1,1,7)$, $(1,4,17)$, $(0,1,20)$, $(1,3,0)$, $(1,20,11)$, $(1,16,18)$, 
$(0,1,21)$, $(1,16,3)$, $(1,4,19)$, $(1,18,11)$, $(0,1,23)$, $(1,19,8)$, $(1,16,8)$, 
$(1,10,23)$, $(1,0,2)$, $(1,19,14)$, $(1,4,1)$, $(1,17,24)$, $(1,0,3)$, $(1,13,20)$, 
$(1,7,23)$, $(1,2,17)$, $(1,0,4)$, $(1,1,17)$, $(1,15,10)$, $(1,15,19)$, $(1,0,6)$, 
$(1,3,13)$, $(1,2,3)$, $(1,5,21)$, $(1,0,8)$, $(1,21,11)$, $(1,6,24)$, $(1,4,22)$, 
$(1,0,10)$, $(1,4,21)$, $(1,13,17)$, $(1,7,1)$, $(1,0,11)$, $(1,15,22)$, $(1,3,2)$, 
$(1,21,18)$, $(1,0,14)$, $(1,17,18)$, $(1,18,12)$, $(1,16,9)$, $(1,0,15)$, $(1,14,8)$, 
$(1,14,16)$, $(1,1,2)$, $(1,0,16)$, $(1,24,15)$, $(1,16,14)$, $(1,6,11)$, $(1,0,17)$, 
$(1,10,16)$, $(1,19,11)$, $(1,24,13)$, $(1,0,18)$, $(1,22,24)$, $(1,5,20)$, $(1,9,6)$, 
$(1,0,21)$, $(1,11,4)$, $(1,11,19)$, $(1,14,15)$, $(1,0,23)$, $(1,16,5)$, $(1,9,21)$, 
$(1,10,5)$, $(1,0,24)$, $(1,18,1)$, $(1,12,18)$, $(1,12,10)$, $(1,1,0)$, $(1,20,0)$, 
$(1,20,13)$, $(1,15,21)$, $(1,1,1)$, $(1,3,19)$, $(1,13,16)$, $(1,7,20)$, $(1,1,3)$, 
$(1,22,1)$, $(1,24,21)$, $(1,18,19)$, $(1,1,8)$, $(1,24,2)$, $(1,5,1)$, $(1,2,6)$, 
$(1,1,9)$, $(1,12,21)$, $(1,9,14)$, $(1,16,12)$, $(1,1,11)$, $(1,9,7)$, $(1,17,23)$, 
$(1,17,3)$, $(1,1,12)$, $(1,8,9)$, $(1,2,11)$, $(1,8,11)$, $(1,1,14)$, $(1,5,5)$, 
$(1,16,6)$, $(1,19,5)$, $(1,1,15)$, $(1,1,18)$, $(1,10,20)$, $(1,22,17)$, $(1,1,21)$, 
$(1,21,3)$, $(1,19,2)$, $(1,3,22)$, $(1,2,2)$, $(1,2,13)$, $(1,22,15)$, $(1,12,19)$, 
$(1,2,4)$, $(1,4,23)$, $(1,11,13)$, $(1,16,4)$, $(1,2,5)$, $(1,17,22)$, $(1,15,5)$, 
$(1,21,20)$, $(1,2,9)$, $(1,18,2)$, $(1,6,10)$, $(1,6,12)$, $(1,2,12)$, $(1,13,9)$, 
$(1,7,16)$, $(1,4,24)$, $(1,2,14)$, $(1,21,10)$, $(1,14,1)$, $(1,10,22)$, $(1,2,18)$, 
$(1,19,7)$, $(1,2,19)$, $(1,10,19)$, $(1,4,0)$, $(1,20,14)$, $(1,2,20)$, $(1,24,0)$, 
$(1,20,8)$, $(1,17,10)$, $(1,2,22)$, $(1,11,24)$, $(1,21,14)$, $(1,24,3)$, $(1,2,23)$, 
$(1,15,12)$, $(1,5,9)$, $(1,22,6)$, $(1,3,1)$, $(1,6,4)$, $(1,9,24)$, $(1,22,19)$, 
$(1,3,4)$, $(1,15,8)$, $(1,18,17)$, $(1,4,11)$, $(1,3,6)$, $(1,16,10)$, $(1,17,15)$, 
$(1,5,13)$, $(1,3,10)$, $(1,8,13)$, $(1,4,14)$, $(1,8,23)$, $(1,3,11)$, $(1,19,1)$, 
$(1,22,0)$, $(1,20,4)$, $(1,3,12)$, $(1,24,9)$, $(1,3,14)$, $(1,18,24)$, $(1,24,4)$, 
$(1,14,20)$, $(1,3,15)$, $(1,10,0)$, $(1,20,1)$, $(1,18,7)$, $(1,3,17)$, $(1,5,22)$, 
$(1,10,6)$, $(1,21,24)$, $(1,3,21)$, $(1,4,12)$, $(1,14,9)$, $(1,9,18)$, $(1,3,23)$, 
$(1,22,20)$, $(1,15,16)$, $(1,11,10)$, $(1,4,5)$, $(1,15,15)$, $(1,6,5)$, $(1,12,8)$, 
$(1,4,6)$, $(1,5,2)$, $(1,4,8)$, $(1,22,18)$, $(1,14,0)$, $(1,20,23)$, $(1,4,9)$, 
$(1,19,3)$, $(1,18,16)$, $(1,22,12)$, $(1,4,20)$, $(1,6,19)$, $(1,22,7)$, $(1,11,11)$, 
$(1,5,3)$, $(1,14,6)$, $(1,15,20)$, $(1,12,11)$, $(1,5,4)$, $(1,5,6)$, $(1,13,6)$, 
$(1,7,7)$, $(1,5,7)$, $(1,21,6)$, $(1,18,15)$, $(1,11,22)$, $(1,5,15)$, $(1,17,6)$, 
$(1,22,14)$, $(1,10,8)$, $(1,5,16)$, $(1,6,6)$, $(1,14,21)$, $(1,16,1)$, $(1,5,23)$, 
$(1,12,6)$, $(1,6,17)$, $(1,19,23)$, $(1,6,0)$, $(1,20,19)$, $(1,9,20)$, $(1,14,24)$, 
$(1,6,2)$, $(1,24,20)$, $(1,17,4)$, $(1,19,21)$, $(1,6,3)$, $(1,17,5)$, $(1,14,10)$, 
$(1,22,5)$, $(1,6,7)$, $(1,8,0)$, $(1,20,18)$, $(1,8,3)$, $(1,6,13)$, $(1,13,2)$, 
$(1,7,24)$, $(1,10,18)$, $(1,6,14)$, $(1,11,14)$, $(1,19,0)$, $(1,20,17)$, $(1,6,15)$, 
$(1,15,17)$, $(1,13,12)$, $(1,7,9)$, $(1,6,20)$, $(1,21,13)$, $(1,10,13)$, $(1,9,22)$, 
$(1,7,0)$, $(1,20,21)$, $(1,12,23)$, $(1,13,4)$, $(1,7,8)$, $(1,9,13)$, $(1,11,16)$, 
$(1,13,18)$, $(1,7,10)$, $(1,19,16)$, $(1,14,7)$, $(1,13,10)$, $(1,7,11)$, $(1,17,14)$, 
$(1,9,4)$, $(1,13,13)$, $(1,7,15)$, $(1,13,11)$, $(1,7,17)$, $(1,12,20)$, $(1,18,8)$, 
$(1,13,14)$, $(1,7,19)$, $(1,24,5)$, $(1,24,18)$, $(1,13,5)$, $(1,7,21)$, $(1,15,7)$, 
$(1,10,14)$, $(1,13,21)$, $(1,7,22)$, $(1,16,23)$, $(1,21,2)$, $(1,13,19)$, $(1,8,1)$, 
$(1,21,17)$, $(1,8,18)$, $(1,16,22)$, $(1,8,4)$, $(1,12,1)$, $(1,8,20)$, $(1,10,3)$, 
$(1,8,8)$, $(1,19,24)$, $(1,8,10)$, $(1,11,2)$, $(1,8,12)$, $(1,9,9)$, $(1,8,17)$, 
$(1,14,4)$, $(1,8,14)$, $(1,22,16)$, $(1,9,8)$, $(1,15,1)$, $(1,9,19)$, $(1,21,8)$, 
$(1,12,7)$, $(1,19,22)$, $(1,9,23)$, $(1,24,11)$, $(1,10,10)$, $(1,18,18)$, $(1,10,9)$, 
$(1,24,8)$, $(1,11,0)$, $(1,20,12)$, $(1,10,17)$, $(1,16,2)$, $(1,14,3)$, $(1,11,3)$, 
$(1,11,18)$, $(1,15,2)$, $(1,17,16)$, $(1,16,19)$, $(1,12,9)$, $(1,16,0)$, $(1,20,7)$, 
$(1,14,12)$, $(1,12,24)$, $(1,15,11)$, $(1,22,21)$, $(1,19,13)$, $(1,14,5)$, $(1,18,13)$, 
$(1,14,13)$, $(1,17,7)$, $(1,24,7)$, $(1,15,24)$, $(1,15,4)$, $(1,16,15)$, $(1,24,12)$, 
$(1,19,9)$, $(1,16,21)$, $(1,22,22)$, $(1,22,10)$, $(1,21,7)$

\subsection{$m_{9}(2,27)\ge 201$}

group: cyclic; generated by 
\[
\begin{pmatrix}
15 &26 &2 \\
12 &3 &17 \\
25 &10 &3 
\end{pmatrix}
\]

\noindent
$(0,0,1)$, $(1,22,6)$, $(1,21,15)$, $(1,21,25)$, $(1,17,10)$, $(1,24,2)$, $(1,22,9)$, 
$(1,24,13)$, $(0,1,1)$, $(1,6,20)$, $(1,7,7)$, $(1,17,17)$, $(1,1,12)$, $(1,0,2)$, 
$(1,25,9)$, $(1,17,16)$, $(0,1,5)$, $(1,2,17)$, $(1,26,14)$, $(1,3,21)$, $(1,5,20)$, 
$(1,14,2)$, $(1,12,9)$, $(1,9,25)$, $(0,1,6)$, $(1,25,16)$, $(1,1,20)$, $(1,15,0)$, 
$(1,10,13)$, $(1,20,2)$, $(1,16,9)$, $(1,18,8)$, $(0,1,14)$, $(1,14,23)$, $(1,4,12)$, 
$(1,16,22)$, $(1,2,1)$, $(1,15,2)$, $(1,13,9)$, $(1,26,26)$, $(0,1,16)$, $(1,19,26)$, 
$(1,18,23)$, $(1,23,3)$, $(1,15,8)$, $(1,12,2)$, $(1,1,9)$, $(1,7,10)$, $(1,0,6)$, 
$(1,16,7)$, $(1,19,20)$, $(1,13,17)$, $(1,14,22)$, $(1,13,18)$, $(1,13,3)$, $(1,11,21)$, 
$(1,0,10)$, $(1,10,10)$, $(1,23,8)$, $(1,26,18)$, $(1,23,11)$, $(1,23,26)$, $(1,11,10)$, 
$(1,13,6)$, $(1,0,11)$, $(1,24,17)$, $(1,22,25)$, $(1,23,24)$, $(1,21,16)$, $(1,3,22)$, 
$(1,26,23)$, $(1,5,10)$, $(1,0,16)$, $(1,8,24)$, $(1,9,22)$, $(1,14,15)$, $(1,11,7)$, 
$(1,15,6)$, $(1,7,20)$, $(1,18,3)$, $(1,0,20)$, $(1,5,12)$, $(1,24,14)$, $(1,22,26)$, 
$(1,7,25)$, $(1,21,20)$, $(1,10,1)$, $(1,2,20)$, $(1,0,21)$, $(1,4,11)$, $(1,26,4)$, 
$(1,4,8)$, $(1,4,1)$, $(1,12,24)$, $(1,25,14)$, $(1,6,17)$, $(1,0,23)$, $(1,14,26)$, 
$(1,3,3)$, $(1,10,14)$, $(1,24,4)$, $(1,22,23)$, $(1,12,21)$, $(1,21,23)$, $(1,1,1)$, 
$(1,8,15)$, $(1,8,11)$, $(1,10,15)$, $(1,13,15)$, $(1,12,14)$, $(1,2,22)$, $(1,8,12)$, 
$(1,1,11)$, $(1,11,19)$, $(1,9,6)$, $(1,25,12)$, $(1,3,8)$, $(1,7,19)$, $(1,13,14)$, 
$(1,4,0)$, $(1,1,14)$, $(1,23,5)$, $(1,2,7)$, $(1,24,0)$, $(1,22,24)$, $(1,16,10)$, 
$(1,21,12)$, $(1,9,20)$, $(1,1,15)$, $(1,18,13)$, $(1,20,12)$, $(1,11,18)$, $(1,4,6)$, 
$(1,3,23)$, $(1,8,4)$, $(1,2,21)$, $(1,1,23)$, $(1,21,1)$, $(1,13,0)$, $(1,16,4)$, 
$(1,7,0)$, $(1,26,0)$, $(1,24,3)$, $(1,22,7)$, $(1,2,5)$, $(1,25,25)$, $(1,11,25)$, 
$(1,12,6)$, $(1,5,0)$, $(1,15,1)$, $(1,5,4)$, $(1,5,25)$, $(1,2,11)$, $(1,12,10)$, 
$(1,8,17)$, $(1,18,0)$, $(1,19,15)$, $(1,23,21)$, $(1,7,13)$, $(1,20,26)$, $(1,4,14)$, 
$(1,5,15)$, $(1,17,22)$, $(1,23,7)$, $(1,16,25)$, $(1,14,4)$, $(1,17,5)$, $(1,8,10)$, 
$(1,4,26)$, $(1,25,26)$, $(1,24,12)$, $(1,22,16)$, $(1,14,3)$, $(1,21,24)$, $(1,18,26)$, 
$(1,12,5)$, $(1,5,5)$, $(1,19,13)$, $(1,20,8)$, $(1,8,25)$, $(1,26,3)$, $(1,16,15)$, 
$(1,10,3)$, $(1,25,3)$, $(1,7,17)$, $(1,25,23)$, $(1,16,13)$, $(1,20,14)$, $(1,18,24)$, 
$(1,11,17)$, $(1,10,6)$, $(1,14,17)$, $(1,9,12)$, $(1,17,21)$, $(1,10,21)$, $(1,19,0)$, 
$(1,11,16)$, $(1,17,6)$, $(1,18,5)$, $(1,26,21)$, $(1,17,19)$

\subsection{$m_{14}(2,29)\ge 364$}

group: cyclic; generated by 
\[
\begin{pmatrix}
0 &14 &20 \\
27 &11 &28 \\
17 &10 &10 
\end{pmatrix}
\]

\noindent
$(0,1,0)$, $(1,7,9)$, $(1,9,7)$, $(1,18,18)$, $(1,11,0)$, $(1,10,27)$, $(1,4,3)$, 
$(0,1,2)$, $(1,5,7)$, $(1,19,3)$, $(1,25,9)$, $(1,4,12)$, $(1,5,15)$, $(1,7,27)$, 
$(0,1,6)$, $(1,18,20)$, $(1,25,18)$, $(1,12,5)$, $(1,22,6)$, $(1,8,28)$, $(1,0,0)$, 
$(0,1,10)$, $(1,8,10)$, $(1,14,5)$, $(1,10,20)$, $(1,21,16)$, $(1,2,2)$, $(1,25,26)$, 
$(0,1,14)$, $(1,21,23)$, $(0,1,17)$, $(1,28,1)$, $(1,17,27)$, $(1,9,13)$, $(1,19,7)$, 
$(0,1,16)$, $(1,19,21)$, $(1,22,25)$, $(1,26,16)$, $(1,26,24)$, $(1,19,8)$, $(1,23,10)$, 
$(0,1,18)$, $(1,12,14)$, $(1,0,28)$, $(1,16,4)$, $(1,8,1)$, $(1,18,23)$, $(1,9,14)$, 
$(0,1,20)$, $(1,17,19)$, $(1,12,0)$, $(1,17,11)$, $(1,0,23)$, $(1,28,18)$, $(1,1,8)$, 
$(0,1,22)$, $(1,4,6)$, $(1,18,15)$, $(1,7,28)$, $(1,27,14)$, $(1,1,17)$, $(1,18,28)$, 
$(0,1,23)$, $(1,22,24)$, $(1,20,20)$, $(1,13,12)$, $(1,28,4)$, $(1,16,24)$, $(1,26,5)$, 
$(0,1,25)$, $(1,9,11)$, $(1,15,22)$, $(1,19,25)$, $(1,20,26)$, $(1,13,11)$, $(1,11,1)$, 
$(0,1,27)$, $(1,14,16)$, $(1,16,10)$, $(1,22,17)$, $(1,2,3)$, $(1,17,9)$, $(1,5,11)$, 
$(1,0,1)$, $(1,10,26)$, $(1,9,0)$, $(1,1,2)$, $(1,20,10)$, $(1,13,18)$, $(1,25,12)$, 
$(1,0,4)$, $(1,5,25)$, $(1,3,6)$, $(1,10,8)$, $(1,10,11)$, $(1,25,2)$, $(1,3,2)$, 
$(1,0,9)$, $(1,3,13)$, $(1,16,22)$, $(1,17,3)$, $(1,1,9)$, $(1,0,16)$, $(1,11,3)$, 
$(1,0,10)$, $(1,4,19)$, $(1,23,15)$, $(1,7,6)$, $(1,22,4)$, $(1,3,12)$, $(1,23,19)$, 
$(1,0,15)$, $(1,7,8)$, $(1,26,12)$, $(1,21,25)$, $(1,8,17)$, $(1,16,14)$, $(1,4,13)$, 
$(1,0,18)$, $(1,8,14)$, $(1,28,10)$, $(1,15,21)$, $(1,4,0)$, $(1,8,15)$, $(1,7,17)$, 
$(1,0,19)$, $(1,22,11)$, $(1,14,24)$, $(1,11,28)$, $(1,26,21)$, $(1,2,23)$, $(1,18,22)$, 
$(1,0,20)$, $(1,23,17)$, $(1,17,21)$, $(1,14,1)$, $(1,25,24)$, $(1,9,4)$, $(1,2,20)$, 
$(1,0,25)$, $(1,21,5)$, $(1,20,18)$, $(1,13,10)$, $(1,3,3)$, $(1,27,9)$, $(1,16,0)$, 
$(1,1,4)$, $(1,9,22)$, $(1,5,16)$, $(1,14,13)$, $(1,8,4)$, $(1,11,6)$, $(1,2,1)$, 
$(1,1,7)$, $(1,26,14)$, $(1,16,26)$, $(1,20,23)$, $(1,13,24)$, $(1,23,27)$, $(1,1,26)$, 
$(1,1,12)$, $(1,17,8)$, $(1,25,21)$, $(1,23,28)$, $(1,21,27)$, $(1,25,16)$, $(1,18,7)$, 
$(1,1,18)$, $(1,3,18)$, $(1,4,23)$, $(1,21,15)$, $(1,7,0)$, $(1,20,0)$, $(1,13,16)$, 
$(1,1,23)$, $(1,16,17)$, $(1,21,20)$, $(1,26,4)$, $(1,15,3)$, $(1,21,9)$, $(1,22,23)$, 
$(1,1,24)$, $(1,14,6)$, $(1,19,5)$, $(1,8,3)$, $(1,2,9)$, $(1,8,8)$, $(1,15,24)$, 
$(1,1,27)$, $(1,23,12)$, $(1,27,7)$, $(1,15,5)$, $(1,5,21)$, $(1,28,14)$, $(1,12,12)$, 
$(1,1,28)$, $(1,8,2)$, $(1,12,25)$, $(1,18,10)$, $(1,28,26)$, $(1,26,25)$, $(1,10,4)$, 
$(1,2,5)$, $(1,23,0)$, $(1,16,5)$, $(1,28,28)$, $(1,14,18)$, $(1,21,21)$, $(1,18,13)$, 
$(1,2,13)$, $(1,11,18)$, $(1,20,16)$, $(1,13,0)$, $(1,14,22)$, $(1,22,0)$, $(1,23,18)$, 
$(1,2,17)$, $(1,19,6)$, $(1,3,20)$, $(1,14,27)$, $(1,14,23)$, $(1,12,7)$, $(1,2,26)$, 
$(1,2,19)$, $(1,15,12)$, $(1,12,23)$, $(1,27,1)$, $(1,14,10)$, $(1,9,12)$, $(1,16,11)$, 
$(1,2,21)$, $(1,20,19)$, $(1,13,4)$, $(1,21,13)$, $(1,14,3)$, $(1,5,9)$, $(1,25,20)$, 
$(1,2,28)$, $(1,28,7)$, $(1,27,28)$, $(1,18,19)$, $(1,14,8)$, $(1,18,26)$, $(1,3,27)$, 
$(1,3,10)$, $(1,5,6)$, $(1,12,17)$, $(1,23,5)$, $(1,19,24)$, $(1,17,13)$, $(1,17,17)$, 
$(1,3,25)$, $(1,9,25)$, $(1,23,21)$, $(1,5,3)$, $(1,16,9)$, $(1,19,26)$, $(1,12,1)$, 
$(1,3,26)$, $(1,11,20)$, $(1,4,22)$, $(1,26,15)$, $(1,7,22)$, $(1,16,21)$, $(1,11,21)$, 
$(1,3,28)$, $(1,12,3)$, $(1,19,9)$, $(1,18,27)$, $(1,22,10)$, $(1,25,7)$, $(1,8,23)$, 
$(1,4,2)$, $(1,27,15)$, $(1,7,16)$, $(1,17,23)$, $(1,23,20)$, $(1,10,25)$, $(1,27,24)$, 
$(1,4,10)$, $(1,11,15)$, $(1,7,10)$, $(1,19,27)$, $(1,28,22)$, $(1,20,7)$, $(1,13,6)$, 
$(1,4,11)$, $(1,4,15)$, $(1,7,15)$, $(1,7,3)$, $(1,10,9)$, $(1,12,4)$, $(1,26,2)$, 
$(1,4,25)$, $(1,15,15)$, $(1,7,5)$, $(1,15,19)$, $(1,26,27)$, $(1,5,5)$, $(1,17,7)$, 
$(1,4,27)$, $(1,10,15)$, $(1,7,7)$, $(1,21,2)$, $(1,19,1)$, $(1,21,11)$, $(1,9,5)$, 
$(1,5,1)$, $(1,27,16)$, $(1,9,6)$, $(1,21,14)$, $(1,27,21)$, $(1,26,8)$, $(1,22,19)$, 
$(1,5,18)$, $(1,22,26)$, $(1,17,22)$, $(1,28,8)$, $(1,8,0)$, $(1,22,12)$, $(1,11,26)$, 
$(1,5,20)$, $(1,20,1)$, $(1,13,14)$, $(1,9,16)$, $(1,10,19)$, $(1,28,6)$, $(1,28,2)$, 
$(1,5,23)$, $(1,11,19)$, $(1,18,24)$, $(1,8,21)$, $(1,21,22)$, $(1,23,11)$, $(1,12,28)$, 
$(1,8,16)$, $(1,25,11)$, $(1,26,0)$, $(1,25,1)$, $(1,15,6)$, $(1,9,18)$, $(1,26,17)$, 
$(1,9,2)$, $(1,28,24)$, $(1,25,22)$, $(1,11,24)$, $(1,27,17)$, $(1,10,1)$, $(1,11,22)$, 
$(1,10,7)$, $(1,16,19)$, $(1,27,2)$, $(1,22,28)$, $(1,10,17)$, $(1,22,27)$, $(1,18,2)$, 
$(1,10,14)$, $(1,15,8)$, $(1,11,2)$, $(1,16,3)$, $(1,15,9)$, $(1,23,22)$, $(1,25,4)$

\subsection{$m_{25}(2,29)\ge 697$}

group: symmetric group $S_3$

\noindent
$(0,0,1)$, $(0,1,0)$, $(1,0,0)$, $(0,1,1)$, $(1,1,0)$, $(1,0,1)$, $(0,1,2)$, 
$(0,1,15)$, $(1,2,0)$, $(1,0,2)$, $(1,0,15)$, $(1,15,0)$, $(0,1,3)$, $(0,1,10)$, 
$(1,3,0)$, $(1,0,3)$, $(1,0,10)$, $(1,10,0)$, $(0,1,4)$, $(0,1,22)$, $(1,4,0)$, 
$(1,0,4)$, $(1,0,22)$, $(1,22,0)$, $(0,1,5)$, $(0,1,6)$, $(1,5,0)$, $(1,0,5)$, 
$(1,0,6)$, $(1,6,0)$, $(0,1,7)$, $(0,1,25)$, $(1,7,0)$, $(1,0,7)$, $(1,0,25)$, 
$(1,25,0)$, $(0,1,9)$, $(0,1,13)$, $(1,9,0)$, $(1,0,9)$, $(1,0,13)$, $(1,13,0)$, 
$(0,1,12)$, $(0,1,17)$, $(1,12,0)$, $(1,0,12)$, $(1,0,17)$, $(1,17,0)$, $(0,1,18)$, 
$(0,1,21)$, $(1,18,0)$, $(1,0,18)$, $(1,0,21)$, $(1,21,0)$, $(0,1,19)$, $(0,1,26)$, 
$(1,19,0)$, $(1,0,19)$, $(1,0,26)$, $(1,26,0)$, $(0,1,23)$, $(0,1,24)$, $(1,23,0)$, 
$(1,0,23)$, $(1,0,24)$, $(1,24,0)$, $(1,1,1)$, $(1,1,2)$, $(1,2,1)$, $(1,15,15)$, 
$(1,1,3)$, $(1,3,1)$, $(1,10,10)$, $(1,1,4)$, $(1,4,1)$, $(1,22,22)$, $(1,1,7)$, 
$(1,7,1)$, $(1,25,25)$, $(1,1,8)$, $(1,8,1)$, $(1,11,11)$, $(1,1,10)$, $(1,10,1)$, 
$(1,3,3)$, $(1,1,11)$, $(1,11,1)$, $(1,8,8)$, $(1,1,12)$, $(1,12,1)$, $(1,17,17)$, 
$(1,1,13)$, $(1,13,1)$, $(1,9,9)$, $(1,1,14)$, $(1,14,1)$, $(1,27,27)$, $(1,1,15)$, 
$(1,15,1)$, $(1,2,2)$, $(1,1,16)$, $(1,16,1)$, $(1,20,20)$, $(1,1,17)$, $(1,17,1)$, 
$(1,12,12)$, $(1,1,18)$, $(1,18,1)$, $(1,21,21)$, $(1,1,19)$, $(1,19,1)$, $(1,26,26)$, 
$(1,1,20)$, $(1,20,1)$, $(1,16,16)$, $(1,1,21)$, $(1,21,1)$, $(1,18,18)$, $(1,1,22)$, 
$(1,22,1)$, $(1,4,4)$, $(1,1,23)$, $(1,23,1)$, $(1,24,24)$, $(1,1,26)$, $(1,26,1)$, 
$(1,19,19)$, $(1,1,27)$, $(1,27,1)$, $(1,14,14)$, $(1,1,28)$, $(1,28,1)$, $(1,28,28)$, 
$(1,2,4)$, $(1,4,2)$, $(1,2,15)$, $(1,15,2)$, $(1,22,15)$, $(1,15,22)$, $(1,2,5)$, 
$(1,5,2)$, $(1,17,15)$, $(1,15,17)$, $(1,6,12)$, $(1,12,6)$, $(1,2,6)$, $(1,6,2)$, 
$(1,3,15)$, $(1,15,3)$, $(1,5,10)$, $(1,10,5)$, $(1,2,8)$, $(1,8,2)$, $(1,4,15)$, 
$(1,15,4)$, $(1,11,22)$, $(1,22,11)$, $(1,2,9)$, $(1,9,2)$, $(1,19,15)$, $(1,15,19)$, 
$(1,13,26)$, $(1,26,13)$, $(1,2,10)$, $(1,10,2)$, $(1,5,15)$, $(1,15,5)$, $(1,3,6)$, 
$(1,6,3)$, $(1,2,12)$, $(1,12,2)$, $(1,6,15)$, $(1,15,6)$, $(1,17,5)$, $(1,5,17)$, 
$(1,2,13)$, $(1,13,2)$, $(1,21,15)$, $(1,15,21)$, $(1,9,18)$, $(1,18,9)$, $(1,2,14)$, 
$(1,14,2)$, $(1,7,15)$, $(1,15,7)$, $(1,27,25)$, $(1,25,27)$, $(1,2,16)$, $(1,16,2)$, 
$(1,8,15)$, $(1,15,8)$, $(1,20,11)$, $(1,11,20)$, $(1,2,17)$, $(1,17,2)$, $(1,23,15)$, 
$(1,15,23)$, $(1,12,24)$, $(1,24,12)$, $(1,2,18)$, $(1,18,2)$, $(1,9,15)$, $(1,15,9)$, 
$(1,21,13)$, $(1,13,21)$, $(1,2,21)$, $(1,21,2)$, $(1,25,15)$, $(1,15,25)$, $(1,18,7)$, 
$(1,7,18)$, $(1,2,22)$, $(1,22,2)$, $(1,11,15)$, $(1,15,11)$, $(1,4,8)$, $(1,8,4)$, 
$(1,2,23)$, $(1,23,2)$, $(1,26,15)$, $(1,15,26)$, $(1,24,19)$, $(1,19,24)$, $(1,2,24)$, 
$(1,24,2)$, $(1,12,15)$, $(1,15,12)$, $(1,23,17)$, $(1,17,23)$, $(1,2,26)$, $(1,26,2)$, 
$(1,13,15)$, $(1,15,13)$, $(1,19,9)$, $(1,9,19)$, $(1,2,27)$, $(1,27,2)$, $(1,28,15)$, 
$(1,15,28)$, $(1,14,28)$, $(1,28,14)$, $(1,2,28)$, $(1,28,2)$, $(1,14,15)$, $(1,15,14)$, 
$(1,28,27)$, $(1,27,28)$, $(1,3,7)$, $(1,7,3)$, $(1,12,10)$, $(1,10,12)$, $(1,25,17)$, 
$(1,17,25)$, $(1,3,8)$, $(1,8,3)$, $(1,22,10)$, $(1,10,22)$, $(1,11,4)$, $(1,4,11)$, 
$(1,3,9)$, $(1,9,3)$, $(1,3,10)$, $(1,10,3)$, $(1,13,10)$, $(1,10,13)$, $(1,3,12)$, 
$(1,12,3)$, $(1,4,10)$, $(1,10,4)$, $(1,17,22)$, $(1,22,17)$, $(1,3,13)$, $(1,13,3)$, 
$(1,14,10)$, $(1,10,14)$, $(1,9,27)$, $(1,27,9)$, $(1,3,14)$, $(1,14,3)$, $(1,24,10)$, 
$(1,10,24)$, $(1,27,23)$, $(1,23,27)$, $(1,3,17)$, $(1,17,3)$, $(1,25,10)$, $(1,10,25)$, 
$(1,12,7)$, $(1,7,12)$, $(1,3,18)$, $(1,18,3)$, $(1,6,10)$, $(1,10,6)$, $(1,21,5)$, 
$(1,5,21)$, $(1,3,19)$, $(1,19,3)$, $(1,16,10)$, $(1,10,16)$, $(1,26,20)$, $(1,20,26)$, 
$(1,3,20)$, $(1,20,3)$, $(1,26,10)$, $(1,10,26)$, $(1,16,19)$, $(1,19,16)$, $(1,3,21)$, 
$(1,21,3)$, $(1,7,10)$, $(1,10,7)$, $(1,18,25)$, $(1,25,18)$, $(1,3,22)$, $(1,22,3)$, 
$(1,17,10)$, $(1,10,17)$, $(1,4,12)$, $(1,12,4)$, $(1,3,23)$, $(1,23,3)$, $(1,27,10)$, 
$(1,10,27)$, $(1,24,14)$, $(1,14,24)$, $(1,3,24)$, $(1,24,3)$, $(1,8,10)$, $(1,10,8)$, 
$(1,23,11)$, $(1,11,23)$, $(1,3,25)$, $(1,25,3)$, $(1,18,10)$, $(1,10,18)$, $(1,7,21)$, 
$(1,21,7)$, $(1,3,26)$, $(1,26,3)$, $(1,28,10)$, $(1,10,28)$, $(1,19,28)$, $(1,28,19)$, 
$(1,3,28)$, $(1,28,3)$, $(1,19,10)$, $(1,10,19)$, $(1,28,26)$, $(1,26,28)$, $(1,4,5)$, 
$(1,5,4)$, $(1,23,22)$, $(1,22,23)$, $(1,6,24)$, $(1,24,6)$, $(1,4,6)$, $(1,6,4)$, 
$(1,16,22)$, $(1,22,16)$, $(1,5,20)$, $(1,20,5)$, $(1,4,9)$, $(1,9,4)$, $(1,24,22)$, 
$(1,22,24)$, $(1,13,23)$, $(1,23,13)$, $(1,4,16)$, $(1,16,4)$, $(1,4,22)$, $(1,22,4)$, 
$(1,20,22)$, $(1,22,20)$, $(1,4,17)$, $(1,17,4)$, $(1,26,22)$, $(1,22,26)$, $(1,12,19)$, 
$(1,19,12)$, $(1,4,20)$, $(1,20,4)$, $(1,5,22)$, $(1,22,5)$, $(1,16,6)$, $(1,6,16)$, 
$(1,4,21)$, $(1,21,4)$, $(1,27,22)$, $(1,22,27)$, $(1,18,14)$, $(1,14,18)$, $(1,4,23)$, 
$(1,23,4)$, $(1,13,22)$, $(1,22,13)$, $(1,24,9)$, $(1,9,24)$, $(1,4,24)$, $(1,24,4)$, 
$(1,6,22)$, $(1,22,6)$, $(1,23,5)$, $(1,5,23)$, $(1,4,25)$, $(1,25,4)$, $(1,28,22)$, 
$(1,22,28)$, $(1,7,28)$, $(1,28,7)$, $(1,4,26)$, $(1,26,4)$, $(1,21,22)$, $(1,22,21)$, 
$(1,19,18)$, $(1,18,19)$, $(1,4,27)$, $(1,27,4)$, $(1,14,22)$, $(1,22,14)$, $(1,14,27)$, 
$(1,27,14)$, $(1,4,28)$, $(1,28,4)$, $(1,7,22)$, $(1,22,7)$, $(1,28,25)$, $(1,25,28)$, 
$(1,5,6)$, $(1,6,5)$, $(1,7,6)$, $(1,6,7)$, $(1,5,25)$, $(1,25,5)$, $(1,5,7)$, 
$(1,7,5)$, $(1,13,6)$, $(1,6,13)$, $(1,25,9)$, $(1,9,25)$, $(1,5,8)$, $(1,8,5)$, 
$(1,19,6)$, $(1,6,19)$, $(1,11,26)$, $(1,26,11)$, $(1,5,9)$, $(1,9,5)$, $(1,25,6)$, 
$(1,6,25)$, $(1,13,7)$, $(1,7,13)$, $(1,5,11)$, $(1,11,5)$, $(1,8,6)$, $(1,6,8)$, 
$(1,8,11)$, $(1,11,8)$, $(1,5,12)$, $(1,12,5)$, $(1,14,6)$, $(1,6,14)$, $(1,17,27)$, 
$(1,27,17)$, $(1,5,13)$, $(1,13,5)$, $(1,20,6)$, $(1,6,20)$, $(1,9,16)$, $(1,16,9)$, 
$(1,5,14)$, $(1,14,5)$, $(1,26,6)$, $(1,6,26)$, $(1,27,19)$, $(1,19,27)$, $(1,5,16)$, 
$(1,16,5)$, $(1,9,6)$, $(1,6,9)$, $(1,20,13)$, $(1,13,20)$, $(1,5,18)$, $(1,18,5)$, 
$(1,21,6)$, $(1,6,21)$, $(1,21,18)$, $(1,18,21)$, $(1,5,19)$, $(1,19,5)$, $(1,27,6)$, 
$(1,6,27)$, $(1,26,14)$, $(1,14,26)$, $(1,5,26)$, $(1,26,5)$, $(1,11,6)$, $(1,6,11)$, 
$(1,19,8)$, $(1,8,19)$, $(1,5,28)$, $(1,28,5)$, $(1,23,6)$, $(1,6,23)$, $(1,28,24)$, 
$(1,24,28)$, $(1,7,8)$, $(1,8,7)$, $(1,26,25)$, $(1,25,26)$, $(1,11,19)$, $(1,19,11)$, 
$(1,7,11)$, $(1,11,7)$, $(1,14,25)$, $(1,25,14)$, $(1,8,27)$, $(1,27,8)$, $(1,7,16)$, 
$(1,16,7)$, $(1,23,25)$, $(1,25,23)$, $(1,20,24)$, $(1,24,20)$, $(1,7,19)$, $(1,19,7)$, 
$(1,11,25)$, $(1,25,11)$, $(1,26,8)$, $(1,8,26)$, $(1,7,20)$, $(1,20,7)$, $(1,7,25)$, 
$(1,25,7)$, $(1,16,25)$, $(1,25,16)$, $(1,7,23)$, $(1,23,7)$, $(1,24,25)$, $(1,25,24)$, 
$(1,24,23)$, $(1,23,24)$, $(1,7,24)$, $(1,24,7)$, $(1,20,25)$, $(1,25,20)$, $(1,23,16)$, 
$(1,16,23)$, $(1,7,27)$, $(1,27,7)$, $(1,8,25)$, $(1,25,8)$, $(1,14,11)$, $(1,11,14)$, 
$(1,8,9)$, $(1,9,8)$, $(1,12,11)$, $(1,11,12)$, $(1,13,17)$, $(1,17,13)$, $(1,8,12)$, 
$(1,12,8)$, $(1,16,11)$, $(1,11,16)$, $(1,17,20)$, $(1,20,17)$, $(1,8,13)$, $(1,13,8)$, 
$(1,27,11)$, $(1,11,27)$, $(1,9,14)$, $(1,14,9)$, $(1,8,14)$, $(1,14,8)$, $(1,9,11)$, 
$(1,11,9)$, $(1,27,13)$, $(1,13,27)$, $(1,8,17)$, $(1,17,8)$, $(1,13,11)$, $(1,11,13)$, 
$(1,12,9)$, $(1,9,12)$, $(1,8,18)$, $(1,18,8)$, $(1,24,11)$, $(1,11,24)$, $(1,21,23)$, 
$(1,23,21)$, $(1,8,21)$, $(1,21,8)$, $(1,28,11)$, $(1,11,28)$, $(1,18,28)$, $(1,28,18)$, 
$(1,8,23)$, $(1,23,8)$, $(1,21,11)$, $(1,11,21)$, $(1,24,18)$, $(1,18,24)$, $(1,9,13)$, 
$(1,13,9)$, $(1,24,13)$, $(1,13,24)$, $(1,9,23)$, $(1,23,9)$, $(1,9,17)$, $(1,17,9)$, 
$(1,18,13)$, $(1,13,18)$, $(1,12,21)$, $(1,21,12)$, $(1,9,21)$, $(1,21,9)$, $(1,12,13)$, 
$(1,13,12)$, $(1,18,17)$, $(1,17,18)$, $(1,9,26)$, $(1,26,9)$, $(1,19,13)$, $(1,13,19)$, 
$(1,19,26)$, $(1,26,19)$, $(1,9,28)$, $(1,28,9)$, $(1,16,13)$, $(1,13,16)$, $(1,28,20)$, 
$(1,20,28)$, $(1,12,17)$, $(1,17,12)$, $(1,28,17)$, $(1,17,28)$, $(1,12,28)$, $(1,28,12)$, 
$(1,12,18)$, $(1,18,12)$, $(1,16,17)$, $(1,17,16)$, $(1,21,20)$, $(1,20,21)$, $(1,12,20)$, 
$(1,20,12)$, $(1,21,17)$, $(1,17,21)$, $(1,16,18)$, $(1,18,16)$, $(1,12,23)$, $(1,23,12)$, 
$(1,14,17)$, $(1,17,14)$, $(1,24,27)$, $(1,27,24)$, $(1,12,27)$, $(1,27,12)$, $(1,24,17)$, 
$(1,17,24)$, $(1,14,23)$, $(1,23,14)$, $(1,14,16)$, $(1,16,14)$, $(1,26,27)$, $(1,27,26)$, 
$(1,20,19)$, $(1,19,20)$, $(1,14,19)$, $(1,19,14)$, $(1,20,27)$, $(1,27,20)$, $(1,26,16)$, 
$(1,16,26)$, $(1,14,20)$, $(1,20,14)$, $(1,18,27)$, $(1,27,18)$, $(1,16,21)$, $(1,21,16)$, 
$(1,14,21)$, $(1,21,14)$, $(1,16,27)$, $(1,27,16)$, $(1,18,20)$, $(1,20,18)$, $(1,16,20)$, 
$(1,20,16)$, $(1,23,20)$, $(1,20,23)$, $(1,16,24)$, $(1,24,16)$, $(1,18,23)$, $(1,23,18)$, 
$(1,19,21)$, $(1,21,19)$, $(1,24,26)$, $(1,26,24)$

\subsection{$m_{25}(2,31)\ge 734$}
group: cyclic; generated by 
\[
\begin{pmatrix}
26 &3 &13 \\
25 &22 &17 \\
5 &23 &14 
\end{pmatrix}
\]

\noindent
$(0,0,1)$, $(1,18,13)$, $(1,22,12)$, $(0,1,0)$, $(1,28,18)$, $(1,16,11)$, $(0,1,2)$, 
$(1,3,21)$, $(1,24,2)$, $(0,1,5)$, $(1,23,0)$, $(1,2,19)$, $(0,1,6)$, $(1,0,4)$, 
$(1,0,29)$, $(0,1,7)$, $(1,17,28)$, $(1,1,24)$, $(0,1,8)$, $(1,29,3)$, $(1,10,10)$, 
$(0,1,9)$, $(1,26,17)$, $(1,26,23)$, $(0,1,10)$, $(1,11,25)$, $(1,6,30)$, $(0,1,11)$, 
$(1,25,1)$, $(1,9,15)$, $(0,1,13)$, $(1,7,23)$, $(1,21,17)$, $(0,1,14)$, $(1,19,29)$, 
$(1,5,4)$, $(0,1,15)$, $(1,5,22)$, $(1,19,27)$, $(0,1,17)$, $(1,9,24)$, $(1,25,28)$, 
$(0,1,18)$, $(1,4,6)$, $(1,8,20)$, $(0,1,19)$, $(1,2,5)$, $(1,23,7)$, $(0,1,20)$, 
$(1,20,14)$, $(1,15,16)$, $(0,1,21)$, $(1,8,8)$, $(1,4,9)$, $(0,1,22)$, $(1,21,30)$, 
$(1,7,25)$, $(0,1,23)$, $(1,22,15)$, $(1,18,1)$, $(0,1,24)$, $(1,24,16)$, $(1,3,14)$, 
$(0,1,27)$, $(1,6,7)$, $(1,11,5)$, $(0,1,28)$, $(1,12,10)$, $(1,30,3)$, $(0,1,29)$, 
$(1,14,11)$, $(1,27,18)$, $(0,1,30)$, $(1,15,27)$, $(1,20,22)$, $(1,0,1)$, $(1,13,14)$, 
$(1,9,27)$, $(1,0,2)$, $(1,22,6)$, $(1,18,19)$, $(1,0,3)$, $(1,15,26)$, $(1,14,26)$, 
$(1,0,5)$, $(1,7,9)$, $(1,30,29)$, $(1,0,6)$, $(1,20,25)$, $(1,26,5)$, $(1,0,7)$, 
$(1,6,3)$, $(1,4,28)$, $(1,0,9)$, $(1,11,2)$, $(1,5,3)$, $(1,0,10)$, $(1,9,21)$, 
$(1,13,20)$, $(1,0,12)$, $(1,28,11)$, $(1,23,18)$, $(1,0,15)$, $(1,8,15)$, $(1,3,22)$, 
$(1,0,16)$, $(1,12,8)$, $(1,2,16)$, $(1,0,17)$, $(1,27,5)$, $(1,25,30)$, $(1,0,18)$, 
$(1,25,24)$, $(1,27,11)$, $(1,0,19)$, $(1,4,22)$, $(1,6,9)$, $(1,0,20)$, $(1,30,23)$, 
$(1,7,15)$, $(1,0,21)$, $(1,16,1)$, $(1,17,13)$, $(1,0,22)$, $(1,29,17)$, $(1,19,25)$, 
$(1,0,23)$, $(1,5,28)$, $(1,11,8)$, $(1,0,24)$, $(1,21,0)$, $(1,10,2)$, $(1,0,26)$, 
$(1,23,12)$, $(1,28,17)$, $(1,0,27)$, $(1,10,27)$, $(1,21,6)$, $(1,0,28)$, $(1,2,10)$, 
$(1,12,14)$, $(1,1,0)$, $(1,23,17)$, $(1,23,15)$, $(1,1,1)$, $(1,3,1)$, $(1,10,7)$, 
$(1,1,2)$, $(1,15,23)$, $(1,30,5)$, $(1,1,3)$, $(1,6,22)$, $(1,25,21)$, $(1,1,6)$, 
$(1,4,8)$, $(1,19,3)$, $(1,1,7)$, $(1,5,15)$, $(1,14,19)$, $(1,1,8)$, $(1,10,19)$, 
$(1,3,17)$, $(1,1,12)$, $(1,28,21)$, $(1,9,4)$, $(1,1,13)$, $(1,12,2)$, $(1,16,25)$, 
$(1,1,16)$, $(1,2,25)$, $(1,27,27)$, $(1,1,19)$, $(1,11,26)$, $(1,11,10)$, $(1,1,21)$, 
$(1,27,14)$, $(1,2,14)$, $(1,1,22)$, $(1,16,30)$, $(1,12,13)$, $(1,1,23)$, $(1,14,16)$, 
$(1,5,23)$, $(1,1,26)$, $(1,20,27)$, $(1,29,2)$, $(1,1,27)$, $(1,9,12)$, $(1,28,30)$, 
$(1,1,28)$, $(1,19,20)$, $(1,4,20)$, $(1,1,29)$, $(1,25,0)$, $(1,6,26)$, $(1,1,30)$, 
$(1,29,28)$, $(1,20,6)$, $(1,2,2)$, $(1,13,19)$, $(1,7,6)$, $(1,2,3)$, $(1,3,0)$, 
$(1,15,3)$, $(1,2,6)$, $(1,21,28)$, $(1,5,30)$, $(1,2,7)$, $(1,29,6)$, $(1,29,21)$, 
$(1,2,9)$, $(1,25,17)$, $(1,10,1)$, $(1,2,11)$, $(1,19,18)$, $(1,14,15)$, $(1,2,12)$, 
$(1,28,1)$, $(1,24,19)$, $(1,2,13)$, $(1,9,30)$, $(1,16,22)$, $(1,2,17)$, $(1,22,2)$, 
$(1,18,29)$, $(1,2,18)$, $(1,26,22)$, $(1,4,11)$, $(1,2,21)$, $(1,14,24)$, $(1,19,17)$, 
$(1,2,22)$, $(1,5,10)$, $(1,21,24)$, $(1,2,23)$, $(1,16,3)$, $(1,9,13)$, $(1,2,24)$, 
$(1,4,5)$, $(1,26,26)$, $(1,2,27)$, $(1,17,8)$, $(1,17,10)$, $(1,2,28)$, $(1,24,12)$, 
$(1,28,2)$, $(1,2,29)$, $(1,11,9)$, $(1,20,5)$, $(1,3,2)$, $(1,30,22)$, $(1,26,20)$, 
$(1,3,3)$, $(1,17,20)$, $(1,9,17)$, $(1,3,4)$, $(1,20,30)$, $(1,21,10)$, $(1,3,5)$, 
$(1,3,25)$, $(1,19,6)$, $(1,3,6)$, $(1,5,11)$, $(1,25,18)$, $(1,3,7)$, $(1,12,24)$, 
$(1,23,14)$, $(1,3,8)$, $(1,11,0)$, $(1,4,7)$, $(1,3,9)$, $(1,22,16)$, $(1,18,4)$, 
$(1,3,10)$, $(1,23,9)$, $(1,12,23)$, $(1,3,18)$, $(1,13,17)$, $(1,6,11)$, $(1,3,19)$, 
$(1,14,10)$, $(1,29,26)$, $(1,3,20)$, $(1,25,26)$, $(1,5,9)$, $(1,3,26)$, $(1,6,4)$, 
$(1,13,25)$, $(1,3,27)$, $(1,26,19)$, $(1,30,28)$, $(1,3,28)$, $(1,21,23)$, $(1,20,8)$, 
$(1,3,29)$, $(1,27,12)$, $(1,28,24)$, $(1,4,0)$, $(1,25,5)$, $(1,29,14)$, $(1,4,2)$, 
$(1,27,16)$, $(1,17,15)$, $(1,4,3)$, $(1,13,1)$, $(1,5,16)$, $(1,4,4)$, $(1,12,11)$, 
$(1,12,18)$, $(1,4,12)$, $(1,28,6)$, $(1,15,10)$, $(1,4,13)$, $(1,10,0)$, $(1,16,28)$, 
$(1,4,14)$, $(1,24,15)$, $(1,9,26)$, $(1,4,15)$, $(1,23,25)$, $(1,7,21)$, $(1,4,19)$, 
$(1,21,14)$, $(1,14,23)$, $(1,4,21)$, $(1,5,19)$, $(1,13,5)$, $(1,4,23)$, $(1,29,27)$, 
$(1,25,4)$, $(1,4,24)$, $(1,20,24)$, $(1,30,1)$, $(1,4,29)$, $(1,4,30)$, $(1,15,12)$, 
$(1,28,27)$, $(1,5,0)$, $(1,29,12)$, $(1,28,9)$, $(1,5,1)$, $(1,12,22)$, $(1,15,0)$, 
$(1,5,2)$, $(1,9,11)$, $(1,10,18)$, $(1,5,7)$, $(1,24,4)$, $(1,6,20)$, $(1,5,8)$, 
$(1,7,14)$, $(1,27,25)$, $(1,5,13)$, $(1,20,10)$, $(1,16,15)$, $(1,5,14)$, $(1,26,1)$, 
$(1,23,27)$, $(1,5,21)$, $(1,30,26)$, $(1,30,8)$, $(1,5,24)$, $(1,22,7)$, $(1,18,14)$, 
$(1,5,25)$, $(1,10,25)$, $(1,9,3)$, $(1,5,26)$, $(1,16,16)$, $(1,20,13)$, $(1,5,27)$, 
$(1,8,28)$, $(1,8,19)$, $(1,6,5)$, $(1,19,23)$, $(1,9,2)$, $(1,6,8)$, $(1,29,30)$, 
$(1,17,1)$, $(1,6,10)$, $(1,27,10)$, $(1,10,29)$, $(1,6,12)$, $(1,28,20)$, $(1,21,16)$, 
$(1,6,13)$, $(1,14,4)$, $(1,16,5)$, $(1,6,14)$, $(1,6,17)$, $(1,8,6)$, $(1,6,16)$, 
$(1,22,22)$, $(1,18,28)$, $(1,6,21)$, $(1,7,27)$, $(1,15,9)$, $(1,6,23)$, $(1,26,0)$, 
$(1,12,21)$, $(1,6,25)$, $(1,23,1)$, $(1,20,20)$, $(1,6,27)$, $(1,16,24)$, $(1,14,13)$, 
$(1,6,29)$, $(1,12,15)$, $(1,26,27)$, $(1,7,1)$, $(1,26,3)$, $(1,8,16)$, $(1,7,2)$, 
$(1,25,12)$, $(1,28,28)$, $(1,7,3)$, $(1,7,19)$, $(1,19,4)$, $(1,7,5)$, $(1,9,1)$, 
$(1,29,10)$, $(1,7,7)$, $(1,8,10)$, $(1,26,2)$, $(1,7,8)$, $(1,22,8)$, $(1,18,22)$, 
$(1,7,11)$, $(1,14,18)$, $(1,17,9)$, $(1,7,12)$, $(1,28,16)$, $(1,25,20)$, $(1,7,16)$, 
$(1,29,7)$, $(1,9,29)$, $(1,7,17)$, $(1,11,14)$, $(1,12,6)$, $(1,7,22)$, $(1,20,26)$, 
$(1,24,7)$, $(1,7,24)$, $(1,24,21)$, $(1,20,17)$, $(1,8,0)$, $(1,22,23)$, $(1,18,21)$, 
$(1,8,2)$, $(1,19,0)$, $(1,27,24)$, $(1,8,4)$, $(1,23,10)$, $(1,14,30)$, $(1,8,9)$, 
$(1,14,3)$, $(1,23,2)$, $(1,8,12)$, $(1,28,7)$, $(1,30,25)$, $(1,8,13)$, $(1,25,15)$, 
$(1,16,10)$, $(1,8,14)$, $(1,21,5)$, $(1,17,0)$, $(1,8,17)$, $(1,24,28)$, $(1,29,4)$, 
$(1,8,18)$, $(1,27,20)$, $(1,19,11)$, $(1,8,22)$, $(1,9,6)$, $(1,11,29)$, $(1,8,23)$, 
$(1,12,29)$, $(1,20,1)$, $(1,8,24)$, $(1,29,25)$, $(1,24,23)$, $(1,8,26)$, $(1,15,21)$, 
$(1,10,8)$, $(1,8,27)$, $(1,11,11)$, $(1,9,18)$, $(1,8,30)$, $(1,10,24)$, $(1,15,20)$, 
$(1,9,7)$, $(1,14,0)$, $(1,9,8)$, $(1,9,9)$, $(1,23,6)$, $(1,30,24)$, $(1,9,20)$, 
$(1,12,9)$, $(1,21,26)$, $(1,9,25)$, $(1,15,11)$, $(1,26,18)$, $(1,9,28)$, $(1,30,21)$, 
$(1,23,29)$, $(1,10,3)$, $(1,24,30)$, $(1,26,29)$, $(1,10,5)$, $(1,14,12)$, $(1,28,22)$, 
$(1,10,6)$, $(1,13,4)$, $(1,30,15)$, $(1,10,12)$, $(1,28,0)$, $(1,14,9)$, $(1,10,14)$, 
$(1,23,22)$, $(1,11,4)$, $(1,10,15)$, $(1,22,14)$, $(1,18,26)$, $(1,10,16)$, $(1,11,19)$, 
$(1,23,24)$, $(1,10,17)$, $(1,12,27)$, $(1,19,7)$, $(1,10,20)$, $(1,20,29)$, $(1,17,14)$, 
$(1,10,21)$, $(1,17,5)$, $(1,20,19)$, $(1,10,22)$, $(1,30,16)$, $(1,13,28)$, $(1,10,28)$, 
$(1,26,15)$, $(1,24,5)$, $(1,11,1)$, $(1,16,21)$, $(1,19,13)$, $(1,11,3)$, $(1,14,29)$, 
$(1,30,30)$, $(1,11,12)$, $(1,28,4)$, $(1,26,21)$, $(1,11,17)$, $(1,29,0)$, $(1,13,15)$, 
$(1,11,18)$, $(1,21,1)$, $(1,25,11)$, $(1,11,22)$, $(1,13,2)$, $(1,29,20)$, $(1,11,23)$, 
$(1,15,25)$, $(1,17,24)$, $(1,11,27)$, $(1,24,20)$, $(1,27,0)$, $(1,11,28)$, $(1,17,17)$, 
$(1,15,4)$, $(1,11,30)$, $(1,27,8)$, $(1,24,1)$, $(1,12,3)$, $(1,27,21)$, $(1,29,22)$, 
$(1,12,7)$, $(1,13,12)$, $(1,28,25)$, $(1,12,12)$, $(1,28,15)$, $(1,13,8)$, $(1,12,20)$, 
$(1,24,8)$, $(1,14,5)$, $(1,12,26)$, $(1,25,2)$, $(1,17,27)$, $(1,12,30)$, $(1,22,20)$, 
$(1,18,24)$, $(1,13,0)$, $(1,21,29)$, $(1,19,26)$, $(1,13,3)$, $(1,16,23)$, $(1,24,13)$, 
$(1,13,6)$, $(1,25,9)$, $(1,15,24)$, $(1,13,7)$, $(1,27,30)$, $(1,13,23)$, $(1,13,10)$, 
$(1,19,8)$, $(1,21,27)$, $(1,13,13)$, $(1,24,14)$, $(1,16,9)$, $(1,13,18)$, $(1,20,3)$, 
$(1,20,11)$, $(1,13,24)$, $(1,26,4)$, $(1,14,8)$, $(1,13,29)$, $(1,23,19)$, $(1,17,25)$, 
$(1,13,30)$, $(1,14,2)$, $(1,26,14)$, $(1,14,7)$, $(1,20,21)$, $(1,25,6)$, $(1,14,14)$, 
$(1,22,29)$, $(1,18,20)$, $(1,14,17)$, $(1,25,10)$, $(1,20,16)$, $(1,14,28)$, $(1,15,5)$, 
$(1,16,4)$, $(1,21,13)$, $(1,15,6)$, $(1,22,0)$, $(1,18,16)$, $(1,15,7)$, $(1,19,2)$, 
$(1,24,10)$, $(1,15,8)$, $(1,23,20)$, $(1,29,5)$, $(1,15,13)$, $(1,21,11)$, $(1,16,18)$, 
$(1,15,15)$, $(1,27,7)$, $(1,15,19)$, $(1,15,18)$, $(1,29,16)$, $(1,23,11)$, $(1,15,30)$, 
$(1,24,9)$, $(1,19,15)$, $(1,16,6)$, $(1,16,13)$, $(1,23,13)$, $(1,16,17)$, $(1,30,13)$, 
$(1,27,17)$, $(1,16,19)$, $(1,27,13)$, $(1,30,20)$, $(1,16,20)$, $(1,26,13)$, $(1,29,19)$, 
$(1,16,26)$, $(1,29,13)$, $(1,26,16)$, $(1,17,4)$, $(1,19,16)$, $(1,26,9)$, $(1,17,11)$, 
$(1,30,18)$, $(1,24,26)$, $(1,17,12)$, $(1,28,12)$, $(1,28,23)$, $(1,17,16)$, $(1,24,0)$, 
$(1,30,6)$, $(1,17,21)$, $(1,22,25)$, $(1,18,15)$, $(1,17,29)$, $(1,26,6)$, $(1,19,22)$, 
$(1,18,5)$, $(1,30,10)$, $(1,22,1)$, $(1,18,9)$, $(1,24,27)$, $(1,22,21)$, $(1,18,11)$, 
$(1,29,18)$, $(1,22,30)$, $(1,18,12)$, $(1,28,26)$, $(1,22,17)$, $(1,18,18)$, $(1,19,5)$, 
$(1,22,11)$, $(1,18,23)$, $(1,23,4)$, $(1,22,27)$, $(1,18,25)$, $(1,20,28)$, $(1,22,9)$, 
$(1,18,30)$, $(1,25,19)$, $(1,22,19)$, $(1,19,10)$, $(1,25,25)$, $(1,30,0)$, $(1,19,30)$, 
$(1,30,7)$, $(1,25,23)$, $(1,21,4)$, $(1,29,1)$, $(1,30,2)$, $(1,21,7)$, $(1,23,16)$, 
$(1,27,29)$, $(1,21,8)$, $(1,27,6)$, $(1,23,3)$, $(1,21,19)$, $(1,24,29)$, $(1,24,25)$, 
$(1,25,8)$, $(1,26,25)$, $(1,25,14)$, $(1,26,28)$, $(1,27,1)$, $(1,27,26)$

\end{document}